\documentclass{amsart}
\usepackage{amsfonts}

\newtheorem{theorem}{Theorem}
\theoremstyle{plain}

\newtheorem{definition}{Definition}

\newtheorem{lemma}{Lemma}

\newtheorem{remark}{Remark}

\numberwithin{equation}{section}
\input{tcilatex}

\begin{document}
\title[Harmonically Convex Functions]{Hermite-Hadamard type inequalities for
harmonically convex functions via fractional integrals}
\author{\.{I}mdat \.{I}\c{s}can}
\address{Department of Mathematics, Faculty of Arts and Sciences,\\
Giresun University, 28100, Giresun, Turkey.}
\email{imdati@yahoo.com, imdat.iscan@giresun.edu.tr}

\begin{abstract}
In this paper, the author established Hermite-Hadamard's inequalities for \
harmonically convex functions via fractional integrals and obtained some
Hermite-Hadamard type inequalities of these classes of functions.
\end{abstract}

\maketitle

\section{Introduction}

Let $f:I\subseteq \mathbb{R\rightarrow R}$ be a convex function defined on
the interval $I$ of real numbers and $a,b\in I$ with $a<b$. The following
inequality%
\begin{equation}
f\left( \frac{a+b}{2}\right) \leq \frac{1}{b-a}\dint\limits_{a}^{b}f(x)dx%
\leq \frac{f(a)+f(b)}{2}  \label{1-1}
\end{equation}

holds. This double inequality is known in the literature as Hermite-Hadamard
integral inequality for convex functions. Note that some of the classical
inequalities for means can be derived from (\ref{1-1}) for appropriate
particular selections of the mapping $f$. Both inequalities hold in the
reversed direction if $f$ is concave. For some results which generalize,
improve and extend the inequalities (\ref{1-1}) we refer the reader to the
recent papers (see \cite{I13,I13a,I13aa,KOA11,SOD10} ) and references
therein.

In \cite{I13}, Iscan gave definition of harmonically convexity as follows:

\begin{definition}
Let $I\subseteq 
\mathbb{R}
\backslash \left\{ 0\right\} $ be a real interval. A function $%
f:I\rightarrow 
\mathbb{R}
$ is said to be harmonically convex, if \ 
\begin{equation}
f\left( \frac{xy}{tx+(1-t)y}\right) \leq tf(y)+(1-t)f(x)  \label{1-2}
\end{equation}%
for all $x,y\in I$ and $t\in \lbrack 0,1]$. If the inequality in (\ref{1-2})
is reversed, then $f$ is said to be harmonically concave.
\end{definition}

The following result of the Hermite-Hadamard type holds.

\begin{theorem}[\protect\cite{I13}]
\label{1.1} Let $f:I\subseteq 
\mathbb{R}
\backslash \left\{ 0\right\} \rightarrow 
\mathbb{R}
$ be a harmonically convex function and $a,b\in I$ with $a<b.$ If $f\in
L[a,b]$ then the following inequalities hold 
\begin{equation*}
f\left( \frac{2ab}{a+b}\right) \leq \frac{ab}{b-a}\dint\limits_{a}^{b}\frac{%
f(x)}{x^{2}}dx\leq \frac{f(a)+f(b)}{2}.
\end{equation*}
\end{theorem}

\begin{lemma}[\protect\cite{I13}]
\label{1.2} Let $f:I\subseteq 
\mathbb{R}
\backslash \left\{ 0\right\} \rightarrow 
\mathbb{R}
$ be a differentiable function on $I^{\circ }$ and $a,b\in I$ with $a<b$. If 
$f^{\prime }\in L[a,b]$ then 
\begin{eqnarray}
&&\frac{f(a)+f(b)}{2}-\frac{ab}{b-a}\dint\limits_{a}^{b}\frac{f(x)}{x^{2}}dx
\label{1-3} \\
&=&\frac{ab\left( b-a\right) }{2}\dint\limits_{0}^{1}\frac{1-2t}{\left(
tb+(1-t)a\right) ^{2}}f^{\prime }\left( \frac{ab}{tb+(1-t)a}\right) dt. 
\notag
\end{eqnarray}
\end{lemma}

In \cite{I13}, Iscan proved the following results connected with the right
part of (\ref{1-2})

\begin{theorem}
\label{1.3}Let $f:I\subseteq \left( 0,\infty \right) \rightarrow 
\mathbb{R}
$ be a differentiable function on $I^{\circ }$, $a,b\in I$ with $a<b,$ and $%
f^{\prime }\in L[a,b].$ If $\left\vert f^{\prime }\right\vert ^{q}$ is
harmonically convex on $[a,b]$ for $q\geq 1,$ then%
\begin{eqnarray}
&&\left\vert \frac{f(a)+f(b)}{2}-\frac{ab}{b-a}\dint\limits_{a}^{b}\frac{f(x)%
}{x^{2}}dx\right\vert  \label{1-4a} \\
&\leq &\frac{ab\left( b-a\right) }{2}\lambda _{1}^{1-\frac{1}{q}}\left[
\lambda _{2}\left\vert f^{\prime }\left( a\right) \right\vert ^{q}+\lambda
_{3}\left\vert f^{\prime }\left( b\right) \right\vert ^{q}\right] ^{\frac{1}{%
q}},  \notag
\end{eqnarray}%
where 
\begin{eqnarray*}
\lambda _{1} &=&\frac{1}{ab}-\frac{2}{\left( b-a\right) ^{2}}\ln \left( 
\frac{\left( a+b\right) ^{2}}{4ab}\right) , \\
\lambda _{2} &=&\frac{-1}{b\left( b-a\right) }+\frac{3a+b}{\left( b-a\right)
^{3}}\ln \left( \frac{\left( a+b\right) ^{2}}{4ab}\right) , \\
\lambda _{3} &=&\frac{1}{a\left( b-a\right) }-\frac{3b+a}{\left( b-a\right)
^{3}}\ln \left( \frac{\left( a+b\right) ^{2}}{4ab}\right) \\
&=&\lambda _{1}-\lambda _{2}.
\end{eqnarray*}
\end{theorem}

\begin{theorem}
\label{1.4}Let $f:I\subseteq \left( 0,\infty \right) \rightarrow 
\mathbb{R}
$ be a differentiable function on $I^{\circ }$, $a,b\in I$ with $a<b,$ and $%
f^{\prime }\in L[a,b].$ If $\left\vert f^{\prime }\right\vert ^{q}$ is
harmonically convex on $[a,b]$ for $q>1,\;\frac{1}{p}+\frac{1}{q}=1,$ then%
\begin{eqnarray}
&&\left\vert \frac{f(a)+f(b)}{2}-\frac{ab}{b-a}\dint\limits_{a}^{b}\frac{f(x)%
}{x^{2}}dx\right\vert  \label{1-4} \\
&\leq &\frac{ab\left( b-a\right) }{2}\left( \frac{1}{p+1}\right) ^{\frac{1}{p%
}}\left( \mu _{1}\left\vert f^{\prime }\left( a\right) \right\vert ^{q}+\mu
_{2}\left\vert f^{\prime }\left( b\right) \right\vert ^{q}\right) ^{\frac{1}{%
q}},  \notag
\end{eqnarray}%
where%
\begin{eqnarray*}
\mu _{1} &=&\frac{\left[ a^{2-2q}+b^{1-2q}\left[ \left( b-a\right) \left(
1-2q\right) -a\right] \right] }{2\left( b-a\right) ^{2}\left( 1-q\right)
\left( 1-2q\right) }, \\
\mu _{2} &=&\frac{\left[ b^{2-2q}-a^{1-2q}\left[ \left( b-a\right) \left(
1-2q\right) +b\right] \right] }{2\left( b-a\right) ^{2}\left( 1-q\right)
\left( 1-2q\right) }.
\end{eqnarray*}
\end{theorem}

We recall the following special functions and inequality

(1) The Beta function:%
\begin{equation*}
\beta \left( x,y\right) =\frac{\Gamma (x)\Gamma (y)}{\Gamma (x+y)}%
=\dint\limits_{0}^{1}t^{x-1}\left( 1-t\right) ^{y-1}dt,\ \ x,y>0,
\end{equation*}

(2) The hypergeometric function:%
\begin{equation*}
_{2}F_{1}\left( a,b;c;z\right) =\frac{1}{\beta \left( b,c-b\right) }%
\dint\limits_{0}^{1}t^{b-1}\left( 1-t\right) ^{c-b-1}\left( 1-zt\right)
^{-a}dt,\ c>b>0,\ \left\vert z\right\vert <1\text{ (see \cite{KST06}).}
\end{equation*}

\begin{lemma}[\protect\cite{PBM81,WZZ13}]
\label{1.5}For $0<\alpha \leq 1$ and $0\leq a<b$, we have%
\begin{equation*}
\left\vert a^{\alpha }-b^{\alpha }\right\vert \leq \left( b-a\right)
^{\alpha }.
\end{equation*}
\end{lemma}

In the following, we will give some necessary definitions and mathematical
preliminaries of fractional calculus theory which are used further in this
paper.

\begin{definition}
Let $f\in L\left[ a,b\right] $. The Riemann-Liouville integrals $%
J_{a+}^{\alpha }f$ and $J_{b-}^{\alpha }f$ of oder $\alpha >0$ with $a\geq 0$
are defined by%
\begin{equation*}
J_{a+}^{\alpha }f(x)=\frac{1}{\Gamma (\alpha )}\dint\limits_{a}^{x}\left(
x-t\right) ^{\alpha -1}f(t)dt,\ x>a
\end{equation*}

and%
\begin{equation*}
J_{b-}^{\alpha }f(x)=\frac{1}{\Gamma (\alpha )}\dint\limits_{x}^{b}\left(
t-x\right) ^{\alpha -1}f(t)dt,\ x<b
\end{equation*}%
respectively, where $\Gamma (\alpha )$ is the Gamma function defined by $%
\Gamma (\alpha )=$ $\dint\limits_{0}^{\infty }e^{-t}t^{\alpha -1}dt$ and $%
J_{a^{+}}^{0}f(x)=J_{b^{-}}^{0}f(x)=f(x).$
\end{definition}

Because of the wide application of Hermite-Hadamard type inequalities and
fractional integrals, many researchers extend their studies to
Hermite-Hadamard type inequalities involving fractional integrals not
limited to integer integrals. Recently, more and more Hermite-Hadamard
inequalities involving fractional integrals have been obtained for different
classes of functions; see \cite{WZZ13,I13b,I13c,I13d,SSYB13,S12,WFZ12}.

The aim of this paper is to establish Hermite--Hadamard's inequalities for
Harmonically convex functions via Riemann--Liouville fractional integral and
some other integral inequalities using the identity is obtained for
fractional integrals.These results have some relationships with \cite{I13}.

\section{\protect\bigskip Main results}

Let $f:I\subseteq \left( 0,\infty \right) \rightarrow 
\mathbb{R}
$ be a differentiable function on $I^{\circ }$, the interior of $I$,
throughout this section we will take%
\begin{eqnarray*}
&&I_{f}\left( g;\alpha ,a,b\right) \\
&=&\frac{f(a)+f(b)}{2}-\frac{\Gamma (\alpha +1)}{2}\left( \frac{ab}{b-a}%
\right) ^{\alpha }\left\{ J_{1/a-}^{\alpha }\left( f\circ g\right)
(1/b)+J_{1/b+}^{\alpha }\left( f\circ g\right) (1/a)\right\}
\end{eqnarray*}%
where $a,b\in I$ with $a<b$, $\alpha >0$, $g(x)=1/x$ and $\Gamma $ is Euler
Gamma function.

Hermite--Hadamard's inequalities for Harmonically convex functions can be
represented in fractional integral forms as follows:

\begin{theorem}
\label{2.0}Let $f:I\subseteq \left( 0,\infty \right) \rightarrow 
\mathbb{R}
$ be a function such that $f\in L[a,b]$, where $a,b\in I$ with $a<b$. If $f$
is a harmonically convex function on $[a,b]$, then the following
inequalities for fractional integrals hold:%
\begin{equation}
f\left( \frac{2ab}{a+b}\right) \leq \frac{\Gamma (\alpha +1)}{2}\left( \frac{%
ab}{b-a}\right) ^{\alpha }\left\{ J_{1/a-}^{\alpha }\left( f\circ g\right)
(1/b)+J_{1/b+}^{\alpha }\left( f\circ g\right) (1/a)\right\} \leq \frac{%
f(a)+f(b)}{2}  \label{2-0}
\end{equation}%
with $\alpha >0$.
\end{theorem}

\begin{proof}
Since $f$ is a harmonically convex function on $[a,b]$, we have for all $%
x,y\in \lbrack a,b]$ (with $t=1/2$ in the inequality (\ref{1-2}))%
\begin{equation*}
f\left( \frac{2xy}{x+y}\right) \leq \frac{f(x)+f(y)}{2}.
\end{equation*}%
Choosing $x=\frac{ab}{tb+(1-t)a}$, $y=\frac{ab}{ta+(1-t)b}$, we get%
\begin{equation}
f\left( \frac{2ab}{a+b}\right) \leq \frac{f\left( \frac{ab}{tb+(1-t)a}%
\right) +f\left( \frac{ab}{ta+(1-t)b}\right) }{2}.  \label{2-0a}
\end{equation}%
Multiplying both sides of (\ref{2-0a}) by $t^{\alpha -1}$, then integrating
the resulting inequality with respect to $t$ over $[0,1]$, we obtain%
\begin{eqnarray*}
f\left( \frac{2ab}{a+b}\right) &\leq &\frac{\alpha }{2}\left\{
\dint\limits_{0}^{1}t^{\alpha -1}f\left( \frac{ab}{tb+(1-t)a}\right)
dt+\dint\limits_{0}^{1}t^{\alpha -1}f\left( \frac{ab}{ta+(1-t)b}\right)
dt\right\} \\
&=&\frac{\alpha }{2}\left( \frac{ab}{b-a}\right) ^{\alpha }\left\{
\dint\limits_{1/b}^{1/a}\left( x-\frac{1}{b}\right) ^{\alpha -1}f\left( 
\frac{1}{x}\right) dx+\dint\limits_{1/b}^{1/a}\left( \frac{1}{a}-x\right)
^{\alpha -1}f\left( \frac{1}{x}\right) dx\right\} \\
&=&\frac{\alpha \Gamma (\alpha )}{2}\left( \frac{ab}{b-a}\right) ^{\alpha
}\left\{ J_{1/a-}^{\alpha }\left( f\circ g\right) (1/b)+J_{1/b+}^{\alpha
}\left( f\circ g\right) (1/a)\right\} \\
&=&\frac{\Gamma (\alpha +1)}{2}\left( \frac{ab}{b-a}\right) ^{\alpha
}\left\{ J_{1/a-}^{\alpha }\left( f\circ g\right) (1/b)+J_{1/b+}^{\alpha
}\left( f\circ g\right) (1/a)\right\} ,\ \text{where }g(x)=1/x.
\end{eqnarray*}%
and the first inequality is proved.

For the proof of the second inequality in (\ref{2-0}) we first note that if $%
f$ is a harmonically convex function, then, for $t\in \left[ 0,1\right] $,
it yields%
\begin{equation*}
f\left( \frac{ab}{tb+(1-t)a}\right) \leq tf(a)+(1-t)f(b)
\end{equation*}%
and%
\begin{equation*}
f\left( \frac{ab}{ta+(1-t)b}\right) \leq tf(b)+(1-t)f(a).
\end{equation*}%
By adding these inequalities we have%
\begin{equation}
f\left( \frac{ab}{tb+(1-t)a}\right) +f\left( \frac{ab}{ta+(1-t)b}\right)
\leq f(a)+f(b).  \label{2-0b}
\end{equation}%
Then multiplying both sides of (\ref{2-0b}) by $t^{\alpha -1}$, and
integrating the resulting inequality with respect to $t$ over $\left[ 0,1%
\right] $, we obtain%
\begin{equation*}
\dint\limits_{0}^{1}f\left( \frac{ab}{tb+(1-t)a}\right) t^{\alpha
-1}dt+\dint\limits_{0}^{1}f\left( \frac{ab}{ta+(1-t)b}\right) t^{\alpha
-1}dt\leq \left[ f(a)+f(b)\right] \dint\limits_{0}^{1}t^{\alpha -1}dt
\end{equation*}%
i.e.%
\begin{equation*}
\Gamma (\alpha +1)\left( \frac{ab}{b-a}\right) ^{\alpha }\left\{
J_{1/a-}^{\alpha }\left( f\circ g\right) (1/b)+J_{1/b+}^{\alpha }\left(
f\circ g\right) (1/a)\right\} \leq f(a)+f(b).
\end{equation*}%
The proof is completed.
\end{proof}

\begin{lemma}
\label{2.1}Let $f:I\subseteq \left( 0,\infty \right) \rightarrow 
\mathbb{R}
$ be a differentiable function on $I^{\circ }$ such that $f^{\prime }\in
L[a,b]$, where $a,b\in I$ with $a<b$. Then the following equality for
fractional integrals holds:%
\begin{eqnarray}
&&I_{f}\left( g;\alpha ,a,b\right)  \label{2-1} \\
&=&\frac{ab\left( b-a\right) }{2}\dint\limits_{0}^{1}\frac{\left[ t^{\alpha
}-(1-t)^{\alpha }\right] }{\left( ta+(1-t)b\right) ^{2}}f^{\prime }\left( 
\frac{ab}{ta+(1-t)b}\right) dt.  \notag
\end{eqnarray}
\end{lemma}

\begin{proof}
Let $A_{t}=ta+(1-t)b$. It suffices to note that 
\begin{eqnarray}
I_{f}\left( g;\alpha ,a,b\right) &=&\frac{ab\left( b-a\right) }{2}%
\dint\limits_{0}^{1}\frac{\left[ t^{\alpha }-(1-t)^{\alpha }\right] }{%
A_{t}^{2}}f^{\prime }\left( \frac{ab}{A_{t}}\right) dt  \notag \\
&=&\frac{ab\left( b-a\right) }{2}\dint\limits_{0}^{1}\frac{t^{\alpha }}{%
A_{t}^{2}}f^{\prime }\left( \frac{ab}{A_{t}}\right) dt-\frac{ab\left(
b-a\right) }{2}\dint\limits_{0}^{1}\frac{(1-t)^{\alpha }}{A_{t}^{2}}%
f^{\prime }\left( \frac{ab}{A_{t}}\right) dt  \notag \\
&&I_{1}+I_{2}.  \label{2-1a}
\end{eqnarray}%
By integrating by part, we have 
\begin{eqnarray}
I_{1} &=&\frac{1}{2}\left[ \left. t^{\alpha }f\left( \frac{ab}{A_{t}}\right)
\right\vert _{0}^{1}-\alpha \dint\limits_{0}^{1}t^{\alpha -1}f\left( \frac{ab%
}{A_{t}}\right) dt\right]  \notag \\
&=&\frac{1}{2}\left[ f\left( b\right) -\alpha \left( \frac{ab}{b-a}\right)
^{\alpha }\dint\limits_{1/b}^{1/a}\left( \frac{1}{a}-x\right) ^{\alpha
-1}f\left( \frac{1}{x}\right) dx\right]  \notag \\
&=&\frac{1}{2}\left[ f\left( b\right) -\Gamma (\alpha +1)\left( \frac{ab}{b-a%
}\right) ^{\alpha }J_{1/b+}^{\alpha }\left( f\circ g\right) (1/a)\right]
\label{2-1b}
\end{eqnarray}%
and similarly we get, 
\begin{eqnarray}
I_{2} &=&-\frac{1}{2}\left[ \left. (1-t)^{\alpha }f\left( \frac{ab}{A_{t}}%
\right) \right\vert _{0}^{1}+\alpha \dint\limits_{0}^{1}(1-t)^{\alpha
-1}f\left( \frac{ab}{A_{t}}\right) dt\right]  \notag \\
&=&-\frac{1}{2}\left[ -f\left( a\right) +\alpha \left( \frac{ab}{b-a}\right)
^{\alpha }\dint\limits_{1/b}^{1/a}(x-\frac{1}{b})^{\alpha -1}f\left( \frac{1%
}{x}\right) dx\right]  \notag \\
&=&\frac{1}{2}\left[ f\left( a\right) -\Gamma (\alpha +1)\left( \frac{ab}{b-a%
}\right) ^{\alpha }J_{1/a-}^{\alpha }\left( f\circ g\right) (1/b)\right] .
\label{2-1c}
\end{eqnarray}%
Using (\ref{2-1b}) and (\ref{2-1c}) in (\ref{2-1a}), we get equality (\ref%
{2-1}).
\end{proof}

\begin{remark}
If Lemma \ref{2.1}, we let $\alpha =1$, then equality (\ref{2-1}) becomes
equality (\ref{1-3}) of Lemma \ref{1.2}.
\end{remark}

Using lemma \ref{2.1}, we can obtain the following fractional integral
inequality.

\begin{theorem}
Let $f:I\subseteq \left( 0,\infty \right) \rightarrow 
\mathbb{R}
$ be a differentiable function on $I^{\circ }$ such that $f^{\prime }\in
L[a,b]$, where $a,b\in I^{\circ }$ with $a<b$. If $\left\vert f^{\prime
}\right\vert ^{q}$ is harmonically convex on $\left[ a,b\right] $ for some
fixed $q\geq 1$, then the following inequality for fractional integrals
holds:%
\begin{eqnarray}
&&\left\vert I_{f}\left( g;\alpha ,a,b\right) \right\vert  \label{2-2} \\
&\leq &\frac{ab\left( b-a\right) }{2}C_{1}^{1-1/q}(\alpha ;a,b)\left(
C_{2}(\alpha ;a,b)\left\vert f^{\prime }(b)\right\vert ^{q}+C_{3}(\alpha
;a,b)\left\vert f^{\prime }(a)\right\vert ^{q}\right) ^{1/q},  \notag
\end{eqnarray}%
where%
\begin{eqnarray*}
C_{1}(\alpha ;a,b) &=&\frac{b^{-2}}{\alpha +1}\left[ _{2}F_{1}\left(
2,1;\alpha +2;1-\frac{a}{b}\right) +_{2}F_{1}\left( 2,\alpha +1;\alpha +2;1-%
\frac{a}{b}\right) \right] , \\
C_{2}(\alpha ;a,b) &=&\frac{b^{-2}}{\alpha +2}\left[ \frac{1}{\alpha +1}%
._{2}F_{1}\left( 2,2;\alpha +3;1-\frac{a}{b}\right) +_{2}F_{1}\left(
2,\alpha +2;\alpha +3;1-\frac{a}{b}\right) \right] , \\
C_{3}(\alpha ;a,b) &=&\frac{b^{-2}}{\alpha +1}\left[ _{2}F_{1}\left(
2,1;\alpha +3;1-\frac{a}{b}\right) +\frac{1}{\alpha +1}._{2}F_{1}\left(
2,\alpha +1;\alpha +3;1-\frac{a}{b}\right) \right] .
\end{eqnarray*}
\end{theorem}

\begin{proof}
Let $A_{t}=ta+(1-t)b$. From Lemma..., using the property of the modulus, the
power mean inequality and the harmonically convexity of $\left\vert
f^{\prime }\right\vert ^{q}$, we find%
\begin{eqnarray*}
&&\left\vert I_{f}\left( g;\alpha ,a,b\right) \right\vert \\
&\leq &\frac{ab\left( b-a\right) }{2}\dint\limits_{0}^{1}\frac{\left\vert
(1-t)^{\alpha }-t^{\alpha }\right\vert }{A_{t}^{2}}\left\vert f^{\prime
}\left( \frac{ab}{A_{t}}\right) \right\vert dt \\
&\leq &\frac{ab\left( b-a\right) }{2}\left( \dint\limits_{0}^{1}\frac{%
\left\vert (1-t)^{\alpha }-t^{\alpha }\right\vert }{A_{t}^{2}}dt\right)
^{1-1/q}\left( \dint\limits_{0}^{1}\frac{\left\vert (1-t)^{\alpha
}-t^{\alpha }\right\vert }{A_{t}^{2}}\left\vert f^{\prime }\left( \frac{ab}{%
A_{t}}\right) \right\vert ^{q}dt\right) ^{1/q}
\end{eqnarray*}%
\begin{equation*}
\leq \frac{ab\left( b-a\right) }{2}\left( \dint\limits_{0}^{1}\frac{\left[
1-t)^{\alpha }+t^{\alpha }\right] }{A_{t}^{2}}dt\right) ^{1-1/q}\left(
\dint\limits_{0}^{1}\frac{(\left[ 1-t)^{\alpha }+t^{\alpha }\right] }{%
A_{t}^{2}}\left[ t\left\vert f^{\prime }(b)\right\vert ^{q}+(1-t)\left\vert
f^{\prime }(a)\right\vert ^{q}\right] dt\right) ^{1/q}
\end{equation*}%
\begin{equation}
\leq \frac{ab\left( b-a\right) }{2}C_{1}^{1-1/q}(\alpha ;a,b)\left(
C_{2}(\alpha ;a,b)\left\vert f^{\prime }(b)\right\vert ^{q}+C_{3}(\alpha
;a,b)\left\vert f^{\prime }(a)\right\vert ^{q}\right) ^{1/q}.  \label{2-2a}
\end{equation}%
calculating $C_{1}(\alpha ;a,b)$, $C_{2}(\alpha ;a,b)$ and $C_{3}(\alpha
;a,b)$, we have 
\begin{eqnarray}
C_{1}(\alpha ;a,b) &=&\dint\limits_{0}^{1}\frac{\left[ 1-t)^{\alpha
}+t^{\alpha }\right] }{A_{t}^{2}}dt  \notag \\
&=&\frac{b^{-2}}{\alpha +1}\left[ _{2}F_{1}\left( 2,1;\alpha +2;1-\frac{a}{b}%
\right) +_{2}F_{1}\left( 2,\alpha +1;\alpha +2;1-\frac{a}{b}\right) \right] ,
\label{2-2b}
\end{eqnarray}%
\begin{eqnarray}
C_{2}(\alpha ;a,b) &=&\dint\limits_{0}^{1}\frac{\left[ 1-t)^{\alpha
}+t^{\alpha }\right] }{A_{t}^{2}}tdt  \notag \\
&=&\frac{b^{-2}}{\alpha +2}\left[ \frac{1}{\alpha +1}._{2}F_{1}\left(
2,2;\alpha +3;1-\frac{a}{b}\right) +_{2}F_{1}\left( 2,\alpha +2;\alpha +3;1-%
\frac{a}{b}\right) \right] ,  \label{2-2c}
\end{eqnarray}%
\begin{eqnarray}
C_{3}(\alpha ;a,b) &=&\dint\limits_{0}^{1}\frac{\left[ 1-t)^{\alpha
}+t^{\alpha }\right] }{A_{t}^{2}}(1-t)dt  \notag \\
&=&\frac{b^{-2}}{\alpha +1}\left[ _{2}F_{1}\left( 2,1;\alpha +3;1-\frac{a}{b}%
\right) +\frac{1}{\alpha +1}._{2}F_{1}\left( 2,\alpha +1;\alpha +3;1-\frac{a%
}{b}\right) \right] ,  \label{2-2d}
\end{eqnarray}%
Thus, if we use (\ref{2-2b}), (\ref{2-2c}) and (\ref{2-2d}) in (\ref{2-2a}),
we obtain the inequality of (\ref{2-2}). This completes the proof.
\end{proof}

When $0<\alpha \leq 1$, using Lemma \ref{1.5} and Lemma \ref{2.1} we shall
give another result for harmonically convex functions as follows.

\begin{theorem}
\label{2.3}Let $f:I\subseteq \left( 0,\infty \right) \rightarrow 
\mathbb{R}
$ be a differentiable function on $I^{\circ }$ such that $f^{\prime }\in
L[a,b]$, where $a,b\in I^{\circ }$ with $a<b$. If $\left\vert f^{\prime
}\right\vert ^{q}$ is harmonically convex on $\left[ a,b\right] $ for some
fixed $q\geq 1$, then the following inequality for fractional integrals
holds:%
\begin{eqnarray}
&&\left\vert I_{f}\left( g;\alpha ,a,b\right) \right\vert   \label{2-3} \\
&\leq &\frac{ab\left( b-a\right) }{2}C_{1}^{1-1/q}(\alpha ;a,b)\left(
C_{2}(\alpha ;a,b)\left\vert f^{\prime }(b)\right\vert ^{q}+C_{3}(\alpha
;a,b)\left\vert f^{\prime }(a)\right\vert ^{q}\right) ^{1/q},  \notag
\end{eqnarray}%
where%
\begin{eqnarray*}
&&C_{1}(\alpha ;a,b) \\
&=&\frac{b^{-2}}{\alpha +1}\left[ _{2}F_{1}\left( 2,\alpha +1;\alpha +2;1-%
\frac{a}{b}\right) -_{2}F_{1}\left( 2,1;\alpha +2;1-\frac{a}{b}\right)
\right.  \\
&&\left. +_{2}F_{1}\left( 2,1;\alpha +2;\frac{1}{2}\left( 1-\frac{a}{b}%
\right) \right) \right] ,
\end{eqnarray*}%
\begin{eqnarray*}
&&C_{2}(\alpha ;a,b) \\
&=&\frac{b^{-2}}{\alpha +2}\left[ _{2}F_{1}\left( 2,\alpha +2;\alpha +3;1-%
\frac{a}{b}\right) -\frac{1}{\alpha +1}._{2}F_{1}\left( 2,2;\alpha +3;1-%
\frac{a}{b}\right) \right.  \\
&&\left. +\frac{1}{2\left( \alpha +1\right) }._{2}F_{1}\left( 2,2;\alpha +3;%
\frac{1}{2}\left( 1-\frac{a}{b}\right) \right) \right] ,
\end{eqnarray*}%
\begin{eqnarray*}
&&C_{3}(\alpha ;a,b) \\
&=&\frac{b^{-2}}{\alpha +2}\left[ \frac{1}{\alpha +1}._{2}F_{1}\left(
2,\alpha +1;\alpha +3;1-\frac{a}{b}\right) -_{2}F_{1}\left( 2,1;\alpha +3;1-%
\frac{a}{b}\right) \right.  \\
&&\left. +_{2}F_{1}\left( 2,1;\alpha +3;\frac{1}{2}\left( 1-\frac{a}{b}%
\right) \right) \right] 
\end{eqnarray*}%
and $0<\alpha \leq 1.$
\end{theorem}

\begin{proof}
Let $A_{t}=ta+(1-t)b$. From Lemma \ref{2.1}, using the property of the
modulus, the power mean inequality and the harmonically convexity of $%
\left\vert f^{\prime }\right\vert ^{q}$, we find%
\begin{eqnarray*}
&&\left\vert I_{f}\left( g;\alpha ,a,b\right) \right\vert \\
&\leq &\frac{ab\left( b-a\right) }{2}\dint\limits_{0}^{1}\frac{\left\vert
(1-t)^{\alpha }-t^{\alpha }\right\vert }{A_{t}^{2}}\left\vert f^{\prime
}\left( \frac{ab}{A_{t}}\right) \right\vert dt \\
&\leq &\frac{ab\left( b-a\right) }{2}\left( \dint\limits_{0}^{1}\frac{%
\left\vert (1-t)^{\alpha }-t^{\alpha }\right\vert }{A_{t}^{2}}dt\right)
^{1-1/q}\left( \dint\limits_{0}^{1}\frac{\left\vert (1-t)^{\alpha
}-t^{\alpha }\right\vert }{A_{t}^{2}}\left\vert f^{\prime }\left( \frac{ab}{%
A_{t}}\right) \right\vert ^{q}dt\right) ^{1/q} \\
&\leq &\frac{ab\left( b-a\right) }{2}K_{1}^{1-1/q}\left( \dint\limits_{0}^{1}%
\frac{\left\vert (1-t)^{\alpha }-t^{\alpha }\right\vert }{A_{t}^{2}}\left[
t\left\vert f^{\prime }(b)\right\vert ^{q}+(1-t)\left\vert f^{\prime
}(a)\right\vert ^{q}\right] dt\right) ^{1/q}
\end{eqnarray*}%
\begin{equation}
\leq \frac{ab\left( b-a\right) }{2}K_{1}^{1-1/q}\left( K_{2}\left\vert
f^{\prime }(b)\right\vert ^{q}+K_{3}\left\vert f^{\prime }(a)\right\vert
^{q}\right) ^{1/q},  \label{2-3a}
\end{equation}%
where 
\begin{eqnarray*}
K_{1} &=&\dint\limits_{0}^{1}\frac{\left\vert (1-t)^{\alpha }-t^{\alpha
}\right\vert }{A_{t}^{2}}dt, \\
K_{2} &=&\dint\limits_{0}^{1}\frac{\left\vert (1-t)^{\alpha }-t^{\alpha
}\right\vert }{A_{t}^{2}}tdt, \\
K_{3} &=&\dint\limits_{0}^{1}\frac{\left\vert (1-t)^{\alpha }-t^{\alpha
}\right\vert }{A_{t}^{2}}(1-t)dt.
\end{eqnarray*}%
Calculating $K_{1}$, $K_{2}$ and $K_{3}$, by Lemma \ref{1.5}, we have%
\begin{eqnarray*}
K_{1} &=&\dint\limits_{0}^{1}\frac{\left\vert (1-t)^{\alpha }-t^{\alpha
}\right\vert }{A_{t}^{2}}dt \\
&=&\dint\limits_{0}^{1/2}\frac{(1-t)^{\alpha }-t^{\alpha }}{A_{t}^{2}}%
dt+\dint\limits_{1/2}^{1}\frac{t^{\alpha }-(1-t)^{\alpha }}{A_{t}^{2}}dt \\
&=&\dint\limits_{0}^{1}\frac{t^{\alpha }-(1-t)^{\alpha }}{A_{t}^{2}}%
dt+2\dint\limits_{0}^{1/2}\frac{(1-t)^{\alpha }-t^{\alpha }}{A_{t}^{2}}dt
\end{eqnarray*}%
\begin{eqnarray*}
&\leq &\dint\limits_{0}^{1}t^{\alpha
}A_{t}^{-2}dt-\dint\limits_{0}^{1}(1-t)^{\alpha
}A_{t}^{-2}dt+2\dint\limits_{0}^{1/2}(1-2t)^{\alpha }A_{t}^{-2}dt \\
&=&\dint\limits_{0}^{1}t^{\alpha
}A_{t}^{-2}dt-\dint\limits_{0}^{1}(1-t)^{\alpha
}A_{t}^{-2}dt+\dint\limits_{0}^{1}(1-u)^{\alpha }b^{-2}\left( 1-u\frac{1}{2}%
(1-\frac{a}{b})\right) ^{-2}du
\end{eqnarray*}%
\begin{eqnarray}
&=&\frac{b^{-2}}{\alpha +1}\left[ _{2}F_{1}\left( 2,\alpha +1;\alpha +2;1-%
\frac{a}{b}\right) -_{2}F_{1}\left( 2,1;\alpha +2;1-\frac{a}{b}\right)
\right.  \notag \\
&&\left. +_{2}F_{1}\left( 2,1;\alpha +2;\frac{1}{2}\left( 1-\frac{a}{b}%
\right) \right) \right]  \notag \\
&=&C_{1}(\alpha ;a,b)  \label{2-3b}
\end{eqnarray}%
and similarly we get%
\begin{eqnarray}
K_{2} &=&\dint\limits_{0}^{1}\frac{\left\vert (1-t)^{\alpha }-t^{\alpha
}\right\vert }{A_{t}^{2}}tdt  \notag \\
&\leq &\dint\limits_{0}^{1}t^{\alpha
+1}A_{t}^{-2}dt-\dint\limits_{0}^{1}(1-t)^{\alpha
}tA_{t}^{-2}dt+2\dint\limits_{0}^{1/2}(1-2t)^{\alpha }tA_{t}^{-2}dt  \notag
\\
&=&\frac{b^{-2}}{\alpha +2}\left[ _{2}F_{1}\left( 2,\alpha +2;\alpha +3;1-%
\frac{a}{b}\right) -\frac{1}{\alpha +1}._{2}F_{1}\left( 2,2;\alpha +3;1-%
\frac{a}{b}\right) \right.  \notag \\
&&\left. +\frac{1}{2\left( \alpha +1\right) }._{2}F_{1}\left( 2,2;\alpha +3;%
\frac{1}{2}\left( 1-\frac{a}{b}\right) \right) \right]  \notag \\
&=&C_{2}(\alpha ;a,b)  \label{2-3c}
\end{eqnarray}%
\begin{eqnarray*}
K_{3} &=&\dint\limits_{0}^{1}\frac{\left\vert (1-t)^{\alpha }-t^{\alpha
}\right\vert }{A_{t}^{2}}(1-t)dt \\
&\leq &\dint\limits_{0}^{1}t^{\alpha
}(1-t)A_{t}^{-2}dt-\dint\limits_{0}^{1}(1-t)^{\alpha
+1}A_{t}^{-2}dt+2\dint\limits_{0}^{1/2}(1-2t)^{\alpha }(1-t)A_{t}^{-2}dt
\end{eqnarray*}%
\begin{eqnarray}
&=&\frac{b^{-2}}{\alpha +2}\left[ \frac{1}{\alpha +1}._{2}F_{1}\left(
2,\alpha +1;\alpha +3;1-\frac{a}{b}\right) \right.  \notag \\
&&\left. -_{2}F_{1}\left( 2,1;\alpha +3;1-\frac{a}{b}\right)
+_{2}F_{1}\left( 2,1;\alpha +3;\frac{1}{2}\left( 1-\frac{a}{b}\right)
\right) \right]  \notag \\
&=&C_{3}(\alpha ;a,b).  \label{2-3d}
\end{eqnarray}

Thus, if we use (\ref{2-3b}), (\ref{2-3c}) and (\ref{2-3d}) in (\ref{2-3a}),
we obtain the inequality of (\ref{2-3}). This completes the proof.
\end{proof}

\begin{remark}
If we take $\alpha =1$ in Theorem \ref{2.3}, then inequality (\ref{2-3})
becomes inequality (\ref{1-4a}) of Theorem \ref{1.3}.
\end{remark}

\begin{theorem}
Let $f:I\subseteq \left( 0,\infty \right) \rightarrow 
\mathbb{R}
$ be a differentiable function on $I^{\circ }$ such that $f^{\prime }\in
L[a,b]$, where $a,b\in I^{\circ }$ with $a<b$. If $\left\vert f^{\prime
}\right\vert ^{q}$ is harmonically convex on $\left[ a,b\right] $ for some
fixed $q>1$, then the following inequality for fractional integrals holds:%
\begin{eqnarray}
&&\left\vert I_{f}\left( g;\alpha ,a,b\right) \right\vert  \label{2-4} \\
&\leq &\frac{a\left( b-a\right) }{2b}\left( \frac{1}{\alpha p+1}\right)
^{1/p}\left( \frac{\left\vert f^{\prime }(b)\right\vert ^{q}+\left\vert
f^{\prime }(a)\right\vert ^{q}}{2}\right) ^{1/q}  \notag \\
&&\times \left[ _{2}F_{1}^{1/p}\left( 2p,1;\alpha p+2;1-\frac{a}{b}\right)
+_{2}F_{1}^{1/p}\left( 2p,\alpha p+1;\alpha p+2;1-\frac{a}{b}\right) \right]
,  \notag
\end{eqnarray}%
where $1/p+1/q=1.$
\end{theorem}

\begin{proof}
Let $A_{t}=ta+(1-t)b$. From Lemma \ref{2.1}, using the H\"{o}lder inequality
and the harmonically convexity of $\left\vert f^{\prime }\right\vert ^{q}$,
we find%
\begin{eqnarray*}
&&\left\vert I_{f}\left( g;\alpha ,a,b\right) \right\vert \\
&\leq &\frac{ab\left( b-a\right) }{2}\dint\limits_{0}^{1}\frac{(1-t)^{\alpha
}}{A_{t}^{2}}\left\vert f^{\prime }\left( \frac{ab}{A_{t}}\right)
\right\vert dt+\dint\limits_{0}^{1}\frac{t^{\alpha }}{A_{t}^{2}}\left\vert
f^{\prime }\left( \frac{ab}{A_{t}}\right) \right\vert dt \\
&\leq &\frac{ab\left( b-a\right) }{2}\left\{ \left( \dint\limits_{0}^{1}%
\frac{(1-t)^{\alpha p}}{A_{t}^{2p}}dt\right) ^{1/p}\left(
\dint\limits_{0}^{1}\left\vert f^{\prime }\left( \frac{ab}{A_{t}}\right)
\right\vert ^{q}dt\right) ^{1/q}\right.
\end{eqnarray*}%
\begin{eqnarray}
&&\left. +\left( \dint\limits_{0}^{1}\frac{t^{\alpha p}}{A_{t}^{2p}}%
dt\right) ^{1/p}\left( \dint\limits_{0}^{1}\left\vert f^{\prime }\left( 
\frac{ab}{A_{t}}\right) \right\vert ^{q}dt\right) ^{1/q}\right\}  \notag \\
&\leq &\frac{ab\left( b-a\right) }{2}\left( K_{4}^{1/p}+K_{5}^{1/p}\right)
\left( \dint\limits_{0}^{1}\left[ t\left\vert f^{\prime }(b)\right\vert
^{q}+(1-t)\left\vert f^{\prime }(a)\right\vert ^{q}\right] dt\right) ^{1/q} 
\notag \\
&\leq &\frac{ab\left( b-a\right) }{2}\left( K_{4}^{1/p}+K_{5}^{1/p}\right)
\left( \frac{\left\vert f^{\prime }(b)\right\vert ^{q}+\left\vert f^{\prime
}(a)\right\vert ^{q}}{2}\right) ^{1/q}.  \label{2-4a}
\end{eqnarray}%
Calculating $K_{4}$ and $K_{5}$, we have 
\begin{eqnarray}
K_{4} &=&\dint\limits_{0}^{1}\frac{(1-t)^{\alpha p}}{A_{t}^{2p}}dt  \notag \\
&=&\frac{b^{-2p}}{\alpha p+1}._{2}F_{1}\left( 2p,1;\alpha p+2;1-\frac{a}{b}%
\right) ,  \label{2-4b}
\end{eqnarray}%
\begin{eqnarray}
K_{5} &=&\dint\limits_{0}^{1}\frac{t^{\alpha p}}{A_{t}^{2p}}dt  \notag \\
&=&\frac{b^{-2p}}{\alpha p+1}._{2}F_{1}\left( 2p,\alpha p+1;\alpha p+2;1-%
\frac{a}{b}\right)  \label{2-4c}
\end{eqnarray}

Thus, if we use (\ref{2-4b}) and (\ref{2-4c}) in (\ref{2-4a}), we obtain the
inequality of (\ref{2-4}). This completes the proof.
\end{proof}

\begin{theorem}
Let $f:I\subseteq \left( 0,\infty \right) \rightarrow 
\mathbb{R}
$ be a differentiable function on $I^{\circ }$ such that $f^{\prime }\in
L[a,b]$, where $a,b\in I^{\circ }$ with $a<b$. If $\left\vert f^{\prime
}\right\vert ^{q}$ is harmonically convex on $\left[ a,b\right] $ for some
fixed $q>1$, then the following inequality for fractional integrals holds:%
\begin{eqnarray}
&&\left\vert I_{f}\left( g;\alpha ,a,b\right) \right\vert  \label{2-5} \\
&\leq &\frac{b-a}{2\left( ab\right) ^{1-1/p}}L_{2p-2}^{2-2/p}(a,b)\left( 
\frac{1}{\alpha q+1}\right) ^{1/q}\left( \frac{\left\vert f^{\prime
}(b)\right\vert ^{q}+\left\vert f^{\prime }(a)\right\vert ^{q}}{2}\right)
^{1/q},  \notag
\end{eqnarray}%
where $1/p+1/q=1$ and $L_{2p-2}(a,b)=\left( \frac{b^{2p-1}-a^{2p-1}}{%
(2p-1)(b-a)}\right) ^{1/(2p-2)}$ is $2p-2$-Logarithmic mean.
\end{theorem}

\begin{proof}
Let $A_{t}=ta+(1-t)b$. From Lemma \ref{2.1} and Lemma \ref{1.5}, using the H%
\"{o}lder inequality and the harmonically convexity of $\left\vert f^{\prime
}\right\vert ^{q}$, we find%
\begin{eqnarray*}
&&\left\vert I_{f}\left( g;\alpha ,a,b\right) \right\vert \\
&\leq &\frac{ab\left( b-a\right) }{2}\dint\limits_{0}^{1}\frac{\left\vert
(1-t)^{\alpha }-t^{\alpha }\right\vert }{A_{t}^{2}}\left\vert f^{\prime
}\left( \frac{ab}{A_{t}}\right) \right\vert dt \\
&\leq &\frac{ab\left( b-a\right) }{2}\left( \dint\limits_{0}^{1}\frac{1}{%
A_{t}^{2p}}dt\right) ^{1/p}\left( \dint\limits_{0}^{1}\left\vert
(1-t)^{\alpha }-t^{\alpha }\right\vert ^{q}\left\vert f^{\prime }\left( 
\frac{ab}{A_{t}}\right) \right\vert ^{q}dt\right) ^{1/q} \\
&\leq &\frac{ab\left( b-a\right) }{2}\left( \dint\limits_{0}^{1}\frac{1}{%
A_{t}^{2p}}dt\right) ^{1/p}\left( \dint\limits_{0}^{1}\left\vert
1-2t\right\vert ^{\alpha q}\left[ t\left\vert f^{\prime }(b)\right\vert
^{q}+(1-t)\left\vert f^{\prime }(a)\right\vert ^{q}\right] dt\right) ^{1/q}
\end{eqnarray*}%
\begin{equation}
\leq \frac{ab\left( b-a\right) }{2}K_{6}^{1/p}\left( K_{7}\left\vert
f^{\prime }(b)\right\vert ^{q}+K_{8}\left\vert f^{\prime }(a)\right\vert
^{q}\right) ^{1/q},  \label{2-5a}
\end{equation}%
where%
\begin{eqnarray}
K_{6} &=&\dint\limits_{0}^{1}\frac{1}{A_{t}^{2p}}dt=b^{-2p}\dint%
\limits_{0}^{1}\left( 1-t\left( 1-\frac{a}{b}\right) \right) ^{-2p}dt  \notag
\\
&=&b^{-2p}._{2}F_{1}\left( 2p,1;2;1-\frac{a}{b}\right) =\frac{%
L_{2p-2}^{2p-2}(a,b)}{\left( ab\right) ^{2p-1}},  \label{2-5b}
\end{eqnarray}%
\begin{eqnarray}
K_{7} &=&\dint\limits_{0}^{1}\left\vert 1-2t\right\vert ^{\alpha q}tdt 
\notag \\
&=&\dint\limits_{0}^{1/2}\left( 1-2t\right) ^{\alpha
q}tdt+\dint\limits_{1/2}^{1}\left( 2t-1\right) ^{\alpha q}tdt  \notag \\
&=&\frac{1}{2\left( \alpha q+1\right) },  \label{2-5c}
\end{eqnarray}%
and 
\begin{eqnarray}
K_{8} &=&\dint\limits_{0}^{1}\left\vert 1-2t\right\vert ^{\alpha q}(1-t)dt 
\notag \\
&=&\frac{1}{2\left( \alpha q+1\right) }.  \label{2-5d}
\end{eqnarray}%
Thus, if we use (\ref{2-5b}), (\ref{2-5c}) and (\ref{2-5d}) in (\ref{2-5a}),
we obtain the inequality of (\ref{2-5}). This completes the proof.
\end{proof}

\begin{theorem}
\label{2.6}Let $f:I\subseteq \left( 0,\infty \right) \rightarrow 
\mathbb{R}
$ be a differentiable function on $I^{\circ }$ such that $f^{\prime }\in
L[a,b]$, where $a,b\in I^{\circ }$ with $a<b$. If $\left\vert f^{\prime
}\right\vert ^{q}$ is harmonically convex on $\left[ a,b\right] $ for some
fixed $q>1$, then the following inequality for fractional integrals holds:%
\begin{eqnarray}
&&\left\vert I_{f}\left( g;\alpha ,a,b\right) \right\vert \leq \frac{a\left(
b-a\right) }{2b}\left( \frac{1}{\alpha p+1}\right) ^{1/p}  \label{2-6} \\
&&\times \left( \frac{_{2}F_{1}\left( 2q,2;3;1-\frac{a}{b}\right) \left\vert
f^{\prime }(b)\right\vert ^{q}+_{2}F_{1}\left( 2q,1;3;1-\frac{a}{b}\right)
\left\vert f^{\prime }(a)\right\vert ^{q}}{2}\right) ^{1/q},  \notag
\end{eqnarray}%
where $1/p+1/q=1$.
\end{theorem}

\begin{proof}
Let $A_{t}=ta+(1-t)b$. From Lemma \ref{2.1} and Lemma \ref{1.5}, using the H%
\"{o}lder inequality and the harmonically convexity of $\left\vert f^{\prime
}\right\vert ^{q}$, we find%
\begin{eqnarray*}
&&\left\vert I_{f}\left( g;\alpha ,a,b\right) \right\vert \\
&\leq &\frac{ab\left( b-a\right) }{2}\dint\limits_{0}^{1}\frac{\left\vert
(1-t)^{\alpha }-t^{\alpha }\right\vert }{A_{t}^{2}}\left\vert f^{\prime
}\left( \frac{ab}{A_{t}}\right) \right\vert dt \\
&\leq &\frac{ab\left( b-a\right) }{2}\left( \dint\limits_{0}^{1}\left\vert
(1-t)^{\alpha }-t^{\alpha }\right\vert ^{p}dt\right) ^{1/p}\left(
\dint\limits_{0}^{1}\frac{1}{A_{t}^{2q}}\left\vert f^{\prime }\left( \frac{ab%
}{A_{t}}\right) \right\vert ^{q}dt\right) ^{1/q} \\
&\leq &\frac{ab\left( b-a\right) }{2}\left( \dint\limits_{0}^{1}\left\vert
1-2t\right\vert ^{\alpha p}dt\right) ^{1/p}\left( \dint\limits_{0}^{1}\frac{1%
}{A_{t}^{2q}}\left[ t\left\vert f^{\prime }(b)\right\vert
^{q}+(1-t)\left\vert f^{\prime }(a)\right\vert ^{q}\right] dt\right) ^{1/q}
\end{eqnarray*}%
\begin{equation}
\leq \frac{ab\left( b-a\right) }{2}K_{9}^{1/p}\left( K_{10}\left\vert
f^{\prime }(b)\right\vert ^{q}+K_{11}\left\vert f^{\prime }(a)\right\vert
^{q}\right) ^{1/q},  \label{2-6a}
\end{equation}%
where%
\begin{equation}
K_{9}=\dint\limits_{0}^{1}\left\vert 1-2t\right\vert ^{\alpha p}dt=\frac{1}{%
\alpha p+1}  \label{2-6b}
\end{equation}%
\begin{eqnarray}
K_{10}
&=&\dint\limits_{0}^{1}tA_{t}^{-2q}dt=b^{-2q}\dint\limits_{0}^{1}t\left(
1-t\left( 1-\frac{a}{b}\right) \right) ^{-2q}dt  \notag \\
&=&\frac{1}{2b^{2q}}._{2}F_{1}\left( 2q,2;3;1-\frac{a}{b}\right)
\label{2-6c}
\end{eqnarray}%
and 
\begin{equation}
K_{11}=\dint\limits_{0}^{1}(1-t)A_{t}^{-2q}dt=\frac{1}{2b^{2q}}%
._{2}F_{1}\left( 2q,1;3;1-\frac{a}{b}\right)  \label{2-6d}
\end{equation}%
Thus, if we use (\ref{2-6b}), (\ref{2-6c}) and (\ref{2-6d}) in (\ref{2-6a}),
we obtain the inequality of (\ref{2-6}). This completes the proof.
\end{proof}

\begin{remark}
If we take $\alpha =1$ in Theorem \ref{2.6}, then inequality (\ref{2-6})
becomes inequality (\ref{1-4}) of Theorem \ref{1.4}.
\end{remark}

\end{document}